\documentclass[12pt]{article}                                               

\input{amssym.def}                                                           
\input{amssym.tex}

\begin{document}                                                             
\title{On  cluster $C^*$-algebras}

\author{Igor ~V. ~Nikolaev}


\date{}
 \maketitle


\newtheorem{thm}{Theorem}
\newtheorem{lem}{Lemma}
\newtheorem{dfn}{Definition}
\newtheorem{rmk}{Remark}
\newtheorem{cor}{Corollary}
\newtheorem{cnj}{Conjecture}
\newtheorem{exm}{Example}


\newcommand{\Coh}{\hbox{\bf Coh}}
\newcommand{\Mod}{\hbox{\bf Mod}}
\newcommand{\Tors}{\hbox{\bf Tors}}

\begin{abstract}
We introduce a   $C^*$-algebra ${\Bbb A}(\mathbf{x}, Q)$ attached to the cluster $\mathbf{x}$ 
and  a  quiver $Q$.   If  $Q_T$ is the quiver coming from   a  triangulation $T$ of  the  Riemann 
surface $S$ with a finite number of cusps,  we  prove  that the 
primitive spectrum of ${\Bbb A}(\mathbf{x}, Q_T)$ times ${\Bbb R}$ 
is homeomorphic  to a generic subset of the  Teichm\"uller space of surface $S$. 
We conclude with an analog of the Tomita-Takesaki theory
and the Connes invariant $T({\cal M})$  for the algebra  
${\Bbb A}(\mathbf{x}, Q_T)$.

\vspace{7mm}


\vspace{5mm}

\end{abstract}

\section{Introduction}
Cluster algebras of rank $m$ are a   class of commutative rings 
 introduced by  [Fomin \& Zelevinsky 2002]  \cite{FoZe1}.
 Among these algebras one finds coordinate rings of  important 
 algebraic varieties,  like the Grassmannians and Schubert varieties;
 cluster algebras appear in the Teichm\"uller theory [Fomin, Shapiro \& Thurston 2008]   \cite{FoShaThu1}.  
 Unlike the coordinate   rings,  the set of generators $x_i$  of  cluster algebra
 is usually infinite and defined by  induction from 
a {\it cluster} $\mathbf{x}=(x_1,\dots,x_m)$    and a 
{\it quiver} $Q$,   see [Williams 2014]  \cite{Wil1} for an excellent  survey;
the cluster algebra is denoted by ${\cal A}(\mathbf{x}, Q)$.   Notice  that   
the ${\cal A}(\mathbf{x}, Q)$ has an additive structure of  countable 
(unperforated)   abelian group  with  an order satisfying the Riesz interpolation property;
see Remark \ref{vershik}.  In other words,   the cluster algebra  ${\cal A}(\mathbf{x}, Q)$ 
is a  {\it dimension group}  by the Effros-Handelman-Shen Theorem  
[Effros 1981, Theorem 3.1]  \cite{E}.

The subject  of our note is an  operator algebra ${\Bbb A}(\mathbf{x}, Q)$,
such that 
\linebreak
$K_0({\Bbb A}(\mathbf{x}, Q))\cong  {\cal A}(\mathbf{x}, Q)$;   
here $K_0({\Bbb A}(\mathbf{x}, Q))$ is the dimension group of  ${\Bbb A}(\mathbf{x}, Q)$
and $\cong$ is an isomorphism of the ordered abelian groups  [Blackadar 1986,  Chapter 7]  \cite{BL}.
The ${\Bbb A}(\mathbf{x}, Q)$ is  an Approximately Finite {\it $C^*$-algebra}
($AF$-algebra) given by a  Bratteli diagram  derived explicitly  from the
pair  $(\mathbf{x}, Q)$.  The  $AF$-algebras were introduced and studied  by  [Bratteli 1972]  \cite{Bra1};
we  refer to  ${\Bbb A}(\mathbf{x}, Q)$  as  a  {\it cluster $C^*$-algebra}.

An exact definition of  ${\Bbb A}(\mathbf{x}, Q)$ can be found   in Section 2.4;  
to  give an idea,  recall that the pair  $(\mathbf{x}, Q)$
is called a {\it seed}  and the cluster algebra  ${\cal A}(\mathbf{x}, Q)$
is generated by seeds obtained 
via   {\it mutation}  of   $(\mathbf{x}, Q)$ (and its  mutants)   in
all  {\it directions}  $k$,   where  $1\le k\le m$  [Williams 2014,  p.5]  \cite{Wil1}.
The mutation process can be described  by an oriented regular tree 
$\overrightarrow{{\Bbb T}}_m$; 
the vertices of  $\overrightarrow{{\Bbb T}}_m$ correspond to the seeds
and the outgoing edges  to the mutations in   
directions  $k$.   
The quotient ${\goth B}(\mathbf{x}, Q)$ of 
$\overrightarrow{{\Bbb T}}_m$ 
 by a relation identifying  equivalent  seeds at  the same level of 
  $\overrightarrow{{\Bbb T}}_m$ is a graph with cycles. 
(For a quick example of such a graph, see Figure 3.) 
The   cluster $C^*$-algebra   ${\Bbb A}(\mathbf{x}, Q)$ is an $AF$-algebra
given by  the ${\goth B}(\mathbf{x}, Q)$ regarded as a  Bratteli diagram   [Bratteli 1972]  \cite{Bra1}.

Let $S_{g,n}$ be a Riemann surface  of genus $g\ge 0$ with $n\ge 1$  cusps and  such that $2g-2+n>0$;
denote by   $T_{g,n}\cong {\Bbb R}^{6g-6+2n}$   the (decorated)  Teichm\"uller space of $S_{g,n}$,
i.e.  a collection   of all Riemann surfaces of genus  $g$ with $n$ cusps endowed with the natural topology  
[Penner 1987]  \cite{Pen1}. 
In what follows,  we   focus on  the algebras  ${\Bbb A}(\mathbf{x}, Q_{g,n})$ with  quivers  $Q_{g,n}$
coming  from  an ideal   triangulation of $S_{g,n}$;   the corresponding cluster  algebra  ${\cal  A}(\mathbf{x}, Q_{g,n})$ 
of  rank $m=6g-6+3n$  is  related  to  the {\it Penner coordinates}  in  $T_{g,n}$  [Fomin, Shapiro \& Thurston 2008]   \cite{FoShaThu1}.

\begin{figure}[here]
\begin{picture}(200,90)(-10,0)

\put(206,20){\circle{5}}
\put(180,60){\circle{5}}
\put(153,20){\circle{5}}

\put(178,60){\vector(2,-3){25}}
\put(182,60){\vector(2,-3){24}}

\put(204,18.5){\vector(-1,0){47}}
\put(204,21.5){\vector(-1,0){47}}

\put(156,22){\vector(2,3){24}}
\put(152,22){\vector(2,3){25}}


\put(141,10){$3$}

\put(210,10){$2$}

\put(177,70){$1$}


\end{picture}
\caption{The Markov quiver $Q_{1,1}$.}
\end{figure}
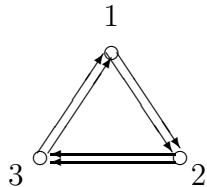

\begin{exm}\label{exm1}
\textnormal{
Let $S_{1,1}$ be a once-punctured torus.    The ideal 
triangulation of $S_{1,1}$ defines   the  Markov quiver
\footnote{Such a quiver is related to solutions in the integer numbers  of the equation $x_1^2+x_2^2+x_3^2=3x_1x_2x_3$ 
considered by A.~A.~Markov;  hence the name.}
$Q_{1,1}$ shown in Figure 1, see  [Fomin, Shapiro \& Thurston 2008,  Example 4.6]   \cite{FoShaThu1}.  
The corresponding cluster $C^*$-algebra ${\Bbb A}(\mathbf{x}, Q_{1,1})$ of rank $3$ can be
written as:   
\begin{equation}\label{eq0}
{\Bbb A}(\mathbf{x}, Q_{1,1})\cong {\goth M}/I_0,
\end{equation}
where $I_0$ is a primitive ideal of an $AF$-algebra  ${\goth M}$.  
The unital $AF$-algebra ${\goth M}$ was originally defined 
by [Mundici 1988,  Section 3]  \cite{Mun1};  the genuine notation 
for such an algebra was  ${\goth M}_1$,  because  $K_0({\goth M}_1)=
(M_1, 1):=$   free one-generator unital  $\ell$-group, 
i.e. a finitely piecewise affine linear continuous real-valued functions on 
$[0,1]$ with integer coefficients.  The   ${\goth M}_1$ was subsequently 
rediscovered after two decades by  [Boca 2008]  \cite{Boc1} and denoted by 
${\goth A}$.   The remarkable properties of ${\goth M}_1$ include the 
following features. Every primitive ideal of ${\goth M}_1$ is essential 
[Mundici 2011, Theorem 4.2]  \cite{Mun3}.  
The ${\goth M}_1$ is equipped with  a faithful invariant tracial state 
[Mundici 2009, Theorem 3.1]  \cite{Mun2}. 
 The center of  ${\goth M}_1$  coincides with the $C^*$-algebra
 $C[0,1]$ of continuous complex valued functions on $[0,1]$  
 [Boca 2008, p. 976]   \cite{Boc1}.  There is an affine weak $\ast$-homeomorphism
 of the state space of $C[0,1]$ onto the space of tracial states on ${\goth M}_1$ 
  [Mundici 2011, Theorem 4.5]  \cite{Mun3}. 
Any state of $C[0,1]$  has precisely one tracial extension to ${\goth M}_1$ 
 [Eckhardt 2011,  Theorem 2.5]  \cite{Eck1}. 
 The automorphism group of ${\goth M}_1$ has precisely two connected 
 components   [Mundici 2011,  Theorem 4.3]  \cite{Mun3}.
 The Gauss map -- a Bernoulli shift for continued fractions -- is generalized 
 in   [Eckhardt 2011]  \cite{Eck1}  to the noncommutative framework of ${\goth M}_1$. 
 In the light of the original definition of ${\goth M}_1$ and the fact that 
 the $K_0$-functor preserves exact sequences (see,  e.g. [Effros 1981, Theorem 3.1] \cite{E}),
 the primitive  spectrum of ${\goth M}_1$ and its hull-kernel topology is widely known 
 to the lattice-ordered group theorists and the MV-algebraists long ago before the 
 laborious analysis in   [Boca 2008]   \cite{Boc1},  where ${\goth M}_1$ is defined 
 in terms of the Bratteli diagram.  We refer the reader to the final part of a paper by 
 [Panti 1999]  \cite{Pan1} for a general result encompassing the characterization of 
 the prime spectrum of $(M_1, 1)\cong Prim~{\goth M}_1$. 
 Moreover,  the  $AF$-algebras  $A_{\theta}$ introduced by   [Effros \& Shen 1980]   \cite{EfSh1}
 are  precisely the infinite-dimensional simple quotients  of ${\goth M}_1$;  this fact 
 was first proved by    [Mundici 1988,  Theorem 3.1(i)]  \cite{Mun1} and rediscovered 
 independently   by [Boca 2008]   \cite{Boc1}. 
 Summing up the above,  the primitive  ideals  $I_{\theta}\subset {\goth M}$ are indexed by  numbers
$\theta\in {\Bbb R}$;  if $\theta$ is irrational,   the quotient ${\goth M}/I_{\theta}\cong A_{\theta}$,  where $A_{\theta}$  
is  the Effros-Shen algebra. 
In view of  (\ref{eq0}),  the algebra ${\goth M}$ is  a non-commutative coordinate
ring of the Teichm\"uller space $T_{1,1}$.   Moreover,  there exists  an analog of   the Tomita-Takesaki
theory  of  modular automorphisms  $\{\sigma_t ~|~ t\in {\Bbb R}\}$  for  algebra   ${\goth M}$,  see Section 4;
such  automorphisms   correspond   to the Teichm\"uller geodesic flow on $T_{1,1}$  [Veech 1986]   \cite{Vee1}. 
The $\sigma_t(I_{\theta})$ is an ideal of ${\goth M}$ for all $t\in {\Bbb R}$,  where $\sigma_0(I_{\theta})=I_{\theta}$. 
The quotient  algebra ${\goth M}/\sigma_t(I_{\theta})$ can be viewed as  a non-commutative 
coordinate ring of the Riemann surface   $S_{1,1}$;  in particular, the pairs $(\theta, t)$  are coordinates in 
the space $T_{1,1}\cong {\Bbb R}^2$.   
We refer the reader to  \cite{Nik2} for  a  construction of the corresponding functor.   
}
\end{exm}

\medskip
Motivated by Example \ref{exm1},  denote by   ${\Bbb A}(\mathbf{x}, Q_{g,n})$
the cluster $C^*$-algebra  corresponding to   a quiver $Q_{g,n}$;
let $\sigma_t:   {\Bbb A}(\mathbf{x}, Q_{g,n})\to {\Bbb A}(\mathbf{x}, Q_{g,n})$
be the Tomita-Takesaki flow on ${\Bbb A}(\mathbf{x}, Q_{g,n})$,  see Section 4
for the details.   Denote by  $Prim~{\Bbb A}(\mathbf{x}, Q_{g,n})$  the set of all primitive ideals of  
${\Bbb A}(\mathbf{x}, Q_{g,n})$ endowed with the Jacobson topology and let $I_{\theta}\in Prim~{\Bbb A}(\mathbf{x}, Q_{g,n})$
for a generic value of index $\theta\in {\Bbb R}^{6g-7+2n}$.  Our main result  can be stated  as follows.    
\begin{thm}\label{thm1}
There exists a homeomorphism
\begin{equation}
h:  Prim~{\Bbb A}(\mathbf{x}, Q_{g,n})\times {\Bbb R}\to \{U\subseteq  T_{g,n} ~|~U~\hbox{{\sf is generic}}\}
\end{equation}
given by the formula $\sigma_t(I_{\theta})\mapsto S_{g,n}$;  the set $U=T_{g,n}$ if and only if
$g=n=1$.   The $\sigma_t(I_{\theta})$
is an ideal of  ${\Bbb A}(\mathbf{x}, Q_{g,n})$ for all $t\in {\Bbb R}$ and 
 the quotient  algebra  ${\Bbb A}(\mathbf{x}, Q_{g,n})/\sigma_t(I_{\theta})$
is  a non-commutative coordinate ring  of  the Riemann surface  $S_{g,n}$.  
\end{thm}
\begin{rmk}
\textnormal{
Theorem \ref{thm1}  is valid for $n\ge 1$,  i.e. the Riemann surfaces with at least one cusp. 
This cannot be improved,  since  the cluster structure of  algebra ${\Bbb A}(\mathbf{x}, Q_{g,n})$
comes from  the Ptolemy relations   satisfied by the Penner coordinates;  so far  such 
coordinates are available  only for the Riemann surfaces with cusps   [Penner 1987]  \cite{Pen1}. 
It is likely,  that the case $n=0$ has also a cluster structure;  
we refer the reader to \cite{Nik1},  where   a  functor from the Riemann surfaces $S_{g,0}$ to
the $AF$-algebras   ${\Bbb A}(\mathbf{x}, Q_{g,0})/\sigma_t(I_{\theta})$  was  constructed. 
}
\end{rmk}
\begin{rmk}
\textnormal{
The braid group $B_{2g+n}$ with $n\in\{1,2\}$ admits a faithful representation by projections 
in  the algebra   ${\Bbb A}(\mathbf{x}, Q_{g,n})$;   such a construction  is based on  the {\it Birman-Hilden Theorem}
for the braid groups.   This observation and the well-known  Laurent phenomenon in  the cluster algebra 
$K_0({\Bbb A}(\mathbf{x}, Q_{g,n}))$ allow  to generalize the Jones and HOMFLY invariants of knots 
and links to an arbitrary number of variables,  see \cite{Nik4} for the details. 
}
\end{rmk}
The article is organized as follows.  We introduce preliminary facts and notation 
in Section 2. Theorem \ref{thm1} is proved in Section 3.  An analog of the 
Tomita-Takesaki theory of modular automorphisms and the Connes invariant $T({\Bbb A}(\mathbf{x}, Q_{g,n}))$  
of  the cluster $C^*$-algebra  ${\Bbb A}(\mathbf{x}, Q_{g,n})$ is  constructed.

\section{Notation}
In this section we introduce notation and briefly review some preliminary facts. 
The reader is encouraged to consult  [Bratteli 1972] \cite{Bra1},  
[Fomin, Shapiro \& Thurston 2008],   [Fomin \& Zelevinsky 2002]  \cite{FoZe1},
[Penner 1987]  \cite{Pen1} and [Williams 2014] \cite{Wil1} for the details.

\subsection{Cluster algebras}
A {\it cluster algebra}  ${\cal A}$ of rank $m$ is a subring of the field ${\Bbb Q}(x_1,\dots, x_m)$ of rational functions
in $n$ variables.  Such an algebra is  defined  by a pair $(\mathbf{x}, B)$,
where $ \mathbf{x}=(x_1,\dots, x_m)$ is a  cluster of  variables and $B=(b_{ij})$ 
is a skew-symmetric integer matrix;   the new cluster $\mathbf{x}'$
is obtained from  $\mathbf{x}$  by an excision of the variable $x_k$ and replacing it by a new variable $x_k'$ 
subject to an exchange relation:
\begin{equation}\label{eq2}
x_kx_k'=\prod_{i=1}^m  x_i^{\max(b_{ik}, 0)} + \prod_{i=1}^m  x_i^{\max(-b_{ik}, 0)}. 
\end{equation}
Since the entries of matrix $B$ are exponents of the monimials in cluster variables, 
one gets  a new pair  $(\mathbf{x}', B')$,  where $B'=(b_{ij}')$ is a skew-symmetric
with: 
\begin{equation}\label{eq3}
b_{ij}'=\cases{
-b_{ij}  & \hbox{if}  $i=k$ \hbox{or} $j=k$\cr
b_{ij}+{|b_{ik}|b_{kj}+b_{ik}|b_{kj}|\over 2}  & \hbox{otherwise.}
}
\end{equation}
For brevity, the pair $(\mathbf{x}, B)$ is called a {\it seed} and the seed $(\mathbf{x}', B'):=(\mathbf{x}',\mu_k(B))$
is obtained from $(\mathbf{x}, B)$ by a {\it mutation} $\mu_k$ in the direction $k$,  where  $1\le k\le m$;  
the $\mu_k$ is an involution, i.e. $\mu_k^2=Id$.   The matrix $B$ is called {\it mutation finite} 
if only finitely many new matrices can be produced from $B$ by repeated matrix mutations. 
 The cluster algebra ${\cal A}(\mathbf{x}, B)$ can be defined
as the subring of ${\Bbb Q}(x_1,\dots, x_m)$ generated by the union of all cluster variables obtained from 
the initial seed $(\mathbf{x}, B)$ by mutations of $(\mathbf{x}, B)$ (and its iterations) in all possible
directions.   We shall write $\overrightarrow{{\Bbb T}}_m$  to denote an oriented tree whose vertices are seeds 
$(\mathbf{x}', B')$ and $m$ outgoing arrows in each vertex correspond to mutations $\mu_k$ of  the seed $(\mathbf{x}', B')$. 
The {\it Laurent phenomenon} 
proved by  [Fomin \& Zelevinsky 2002]  \cite{FoZe1}  says  that  ${\cal A}(\mathbf{x}, B)\subset {\Bbb Z}[\mathbf{x}^{\pm 1}]$,
where  ${\Bbb Z}[\mathbf{x}^{\pm 1}]$ is the ring of  the Laurent polynomials in  variables $\mathbf{x}=(x_1,\dots,x_n)$;  
in other words, each  generator $x_i$  of  algebra ${\cal A}(\mathbf{x}, B)$  can be 
written as a  Laurent polynomial in $n$ variables with the   integer coefficients. 
\begin{rmk}\label{vershik}
\textnormal{
The Laurent phenomenon 
 turns the additive structure of  cluster algebra  ${\cal A}(\mathbf{x}, B)$ into  a totally ordered  abelian group 
satisfying the Riesz interpolation property,  
i.e. a   dimension group   [Effros 1981, Theorem 3.1] \cite{E};
the    abelian group with order comes  from the semigroup of the Laurent polynomials with {\it positive} 
coefficients,     see  \cite{Nik3} for the details. A background on the partially and totally ordered, unperforated
abelian groups with the Riesz interpolation property  can be found in   [Effros 1981]  \cite{E}.  
}
\end{rmk}
To deal with mutation formulas (\ref{eq2}) and (\ref{eq3}) in geometric terms, recall that  
a {\it quiver} $Q$ is an oriented graph given by the set of vertices $Q_0$ and the set
of arrows $Q_1$;  an example of quiver is given in Figure 1.  
Let $k$ be a vertex of $Q$;  the mutated at vertex $k$ quiver $\mu_k(Q)$  has the
same set of vertices as $Q$  but the set of arrows is obtained by the following procedure:
(i) for each sub-quiver $i\to k\to j$ one adds a new arrow $i\to j$; (ii) one reverses all arrows 
with source or target $k$; (iii) one removes the arrows in a maximal set of pairwise disjoint 
$2$-cycles.   The reader can verify,  that  if one encodes a quiver $Q$ with $n$ vertices
by a skew-symmetric matrix $B(Q)=(b_{ij})$ with $b_{ij}$ equal to to the number of arrows 
from vertex $i$ to vertex $j$,  then  mutation $\mu_k$  of seed   $(\mathbf{x}, B)$ coincides with 
such  of the corresponding quiver $Q$. Thus the cluster algebra  ${\cal A}(\mathbf{x}, B)$
is defined by a quiver $Q$;  we shall denote  such an algebra by ${\cal A}(\mathbf{x}, Q)$.

\subsection{Cluster algebras from Riemann surfaces}
Let $g$ and $n$ be integers, such that $g\ge 0, ~n\ge 1$ and $2g-2+n>0$. 
Denote by $S_{g,n}$ a Riemann surface of genus $g$ with the $n$ cusp points.
It is known  that the fundamental domain of $S_{g,n}$ can be triangulated by 
$6g-6+3n$ geodesic arcs $\gamma$,  such that  the footpoints  of  each arc
at the absolute of Lobachevsky plane ${\Bbb H}=\{x+iy\in {\Bbb C} ~|~ y>0\}$
coincide with a (pre-image of) cusp of the $S_{g,n}$.  If $l(\gamma)$ is the
hyperbolic length of $\gamma$ measured (with a sign) between two horocycles
around the footpoints of $\gamma$,  then we set $\lambda(\gamma)=e^{{1\over 2} l(\gamma)}$;
the $\lambda(\gamma)$ are known to satisfy the {\it Ptolemy relation}:     
\begin{equation}\label{eq4}
\lambda(\gamma_1)\lambda(\gamma_2)+\lambda(\gamma_3)\lambda(\gamma_4)=
\lambda(\gamma_5)\lambda(\gamma_6),
\end{equation}
where $\gamma_1, \dots, \gamma_4$ are pairwise opposite sides and $\gamma_5, \gamma_6$
are the diagonals of a geodesic quadrilateral in ${\Bbb H}$.

Denote by $T_{g,n}$ the decorated Teichm\"uller space of $S_{g,n}$, i.e. the set of all 
complex surfaces of genus $g$ with $n$ cusps endowed with the natural topology;
it is known that $T_{g,n}\cong {\Bbb R}^{6g-6+2n}$.  
\begin{thm}\label{thm2}
{\bf ([Penner 1987]  \cite{Pen1})}  
The map $\lambda$ on  the set of $6g-6+3n$ geodesic arcs $\gamma_i$
defining a triangulation of $S_{g,n}$ is  a homeomorphism with the image $T_{g,n}$.   
\end{thm}
\begin{rmk}
\textnormal{
Notice that  among $6g-6+3n$ real numbers $\lambda(\gamma_i)$ there are only $6g-6+2n$ independent, 
since such numbers  must satisfy $n$ Ptolemy relations  (\ref{eq4}). 
}
\end{rmk}

Let $T$ be a triangulation of surface $S_{g,n}$ by $6g-6+3n$ geodesic arcs $\gamma_i$; 
consider a skew-symmetric matrix $B_T=(b_{ij})$,  where  $b_{ij}$ is equal to the number 
of triangles in $T$ with sides $\gamma_i$ and $\gamma_j$ in clockwise order minus 
the number of triangles in $T$ with  sides $\gamma_i$ and $\gamma_j$ in the 
counter-clockwise order.  It is known  that  matrix $B_T$ is always mutation finite. 
The cluster algebra ${\cal A}(\mathbf{x}, B_T)$ 
of rank $6g-6+3n$ is called {\it associated} to triangulation $T$.  
\begin{exm}\label{exm2}
\textnormal{
Let  $S_{1,1}$ be  a once-puncuted torus of Example \ref{exm1}. 
The triangulation $T$ of the fundamental domain $({\Bbb R}^2 - {\Bbb Z}^2)/{\Bbb Z}^2$
of $S_{1,1}$ is sketched  in Figure 2 in the charts ${\Bbb R}^2$ and ${\Bbb H}$, respectively. 
 It is easy to see that in this case $\mathbf{x}=(x_1,x_2,x_3)$ with $x_1=\gamma_{23},
 x_2=\gamma_{34}$ and $x_3=\gamma_{24}$,  where $\gamma_{ij}$ denotes
 a geodesic arc with the footpoints $i$ and $j$.  The Ptolemy relation (\ref{eq4}) 
 reduces to $\lambda^2(\gamma_{23})+\lambda^2(\gamma_{34})=\lambda^2(\gamma_{24})$;
 thus $T_{1,1}\cong {\Bbb R}^2$.   The reader is encouraged to verify, that matrix $B_T$ 
  has the form: 
 \begin{equation}\label{eq5}
 B_T=\left(
 \matrix{0 & 2 & -2\cr
              -2 & 0 & 2\cr
               2 & -2 & 0}
              \right).
 \end{equation}
  }
\end{exm}
\begin{figure}[here]
\begin{picture}(200,100)(0,0)


\put(100,30){\circle{4}}
\put(140,30){\circle{4}}
\put(100,58){\circle{4}}
\put(140,58){\circle{4}}

\put(101,31){\line(3,2){38}}
\put(101,30){\line(1,0){38}}
\put(100,31){\line(0,1){25}}

\put(101,58){\line(1,0){38}}
\put(140,31){\line(0,1){25}}

\put(85,55){$1$}
\put(85,30){$4$}
\put(148,55){$2$}
\put(148,30){$3$}

\put(115,10){${\Bbb R}^2$}


\put(200,30){\line(1,0){75}}
\put(210,30){\line(0,1){40}}
\put(260,30){\line(0,1){40}}

\qbezier(240,30)(250,47)(260,30)
\qbezier(240,30)(225,60)(210,30)
\qbezier(210,30)(235,90)(260,30)

\put(218,75){$1=\infty$}
\put(228,32){$3$}
\put(265,32){$2$}
\put(200,32){$4$}

\put(230,10){${\Bbb H}$}

\end{picture}
\caption{Triangulation of the Riemann surface $S_{1,1}$.}
\end{figure}
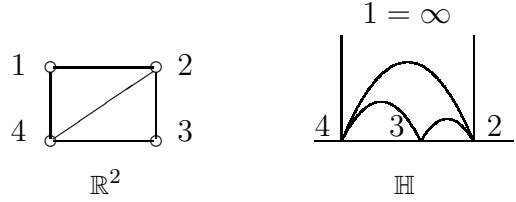
\begin{thm}\label{thm3}
{\bf ([Fomin, Shapiro \& Thurston 2008]  \cite{FoShaThu1})}  
The cluster algebra ${\cal A}(\mathbf{x}, B_T)$ does not depend on triangulation
$T$, but only on the surface $S_{g,n}$;  namely,   replacement of the geodesic 
arc $\gamma_k$ by a new geodesic arc $\gamma_k'$ (a flip of $\gamma_k$) 
corresponds to a mutation $\mu_k$ of the seed $(\mathbf{x}, B_T)$. 
\end{thm}
\begin{rmk}
\textnormal{
In view of Theorems \ref{thm2} and \ref{thm3},  the  ${\cal A}(\mathbf{x}, B_T)$
corresponds to an algebra of functions on the Teichm\"uller space $T_{g,n}$;
such an algebra is an analog of the  coordinate ring of $T_{g,n}$.  
}
\end{rmk}

\subsection{$C^*$-algebras}
A {\it $C^*$-algebra} is an algebra $A$ over ${\Bbb C}$ with a norm
$a\mapsto ||a||$ and an involution $a\mapsto a^*$ such that
it is complete with respect to the norm and $||ab||\le ||a||~ ||b||$
and $||a^*a||=||a^2||$ for all $a,b\in A$.
Any commutative $C^*$-algebra is  isomorphic
to the algebra $C_0(X)$ of continuous complex-valued
functions on some locally compact Hausdorff space $X$; 
otherwise, $A$ represents a noncommutative  topological
space.

An {\it $AF$-algebra}  (Approximately Finite $C^*$-algebra) is defined to
be the  norm closure of an ascending sequence of   finite dimensional
$C^*$-algebras $M_n$,  where  $M_n$ is the $C^*$-algebra of the $n\times n$ matrices
with entries in ${\Bbb C}$. Here the index $n=(n_1,\dots,n_k)$ represents
the  semi-simple matrix algebra $M_n=M_{n_1}\oplus\dots\oplus M_{n_k}$.
The ascending sequence mentioned above  can be written as 
\begin{equation}
M_1\buildrel\rm\varphi_1\over\longrightarrow M_2
   \buildrel\rm\varphi_2\over\longrightarrow\dots,
\end{equation}
where $M_i$ are the finite dimensional $C^*$-algebras and
$\varphi_i$ the homomorphisms between such algebras.  
The homomorphisms $\varphi_i$ can be arranged into  a graph as follows. 
Let  $M_i=M_{i_1}\oplus\dots\oplus M_{i_k}$ and 
$M_{i'}=M_{i_1'}\oplus\dots\oplus M_{i_k'}$ be 
the semi-simple $C^*$-algebras and $\varphi_i: M_i\to M_{i'}$ the  homomorphism. 
One has  two sets of vertices $V_{i_1},\dots, V_{i_k}$ and $V_{i_1'},\dots, V_{i_k'}$
joined by  $b_{rs}$ edges  whenever the summand $M_{i_r}$ contains $b_{rs}$
copies of the summand $M_{i_s'}$ under the embedding $\varphi_i$. 
As $i$ varies, one obtains an infinite graph called the  {\it Bratteli diagram} of the
$AF$-algebra.  The matrix $B=(b_{rs})$ is known as  a {\it partial multiplicity} matrix;
an infinite sequence of $B_i$ defines a unique $AF$-algebra.

\medskip
Let $\theta\in {\Bbb R}^{n-1}$;  recall that by the {\it Jacobi-Perron continued fraction} of 
vector $(1,\theta)$ one understands the limit:
$$
\left(\matrix{1\cr \theta_1\cr\vdots\cr\theta_{n-1}} \right)=
\lim_{k\to\infty} 
\left(\matrix{0 &  0 & \dots & 0 & 1\cr
              1 &  0 & \dots & 0 & b_1^{(1)}\cr
              \vdots &\vdots & &\vdots &\vdots\cr
              0 &  0 & \dots & 1 & b_{n-1}^{(1)}}\right)
\dots 
\left(\matrix{0 &  0 & \dots & 0 & 1\cr
              1 &  0 & \dots & 0 & b_1^{(k)}\cr
              \vdots &\vdots & &\vdots &\vdots\cr
              0 &  0 & \dots & 1 & b_{n-1}^{(k)}}\right)
\left(\matrix{0\cr 0\cr\vdots\cr 1} \right),
$$
where $b_i^{(j)}\in {\Bbb N}\cup\{0\}$, see e.g.  [Bernstein 1971] \cite{B};
the limit converges for a generic subset of vectors $\theta\in {\Bbb R}^{n-1}$.
Notice that $n=2$ corresponds to (a matrix form of) the regular continued
fraction of $\theta$;  such a fraction is always convergent. 
Moreover,  the Jacobi-Perron fraction is {\it finite} if and only if 
vector $\theta=(\theta_i)$,  where $\theta_i$ are {\it rational}.   
The $AF$-algebra  $A_{\theta}$  {\it associated} to the vector $(1,\theta)$ 
is defined by the Bratteli diagram with the partial multiplicity matrices equal to  $B_k$ 
 in the Jacobi-Perron fraction of $(1,\theta)$;  in particular, if $n=2$ the $A_{\theta}$ 
 coincides with the Effros-Shen algebra  [Effros \& Shen 1980]   \cite{EfSh1}.

\subsection{Cluster $C^*$-algebras}
Notice that the mutation tree $\overrightarrow{{\Bbb T}}_m$  of 
a cluster algebra ${\cal A}(\mathbf{x}, B)$ has a grading by  levels,
i.e. a distance from the root of $\overrightarrow{{\Bbb T}}_m$. 
We shall say that a pair of  clusters $\mathbf{x}$ and $\mathbf{x}'$ are  
{\it $\ell$-equivalent},    if:

\medskip
  (i)  $\mathbf{x}$ and $\mathbf{x}'$ lie at the same level;   

\smallskip
(ii) $\mathbf{x}$ and $\mathbf{x}'$  coincide
modulo a cyclic permutation of variables $x_i$;

\smallskip
(iii)  $B=B'$.  

\medskip\noindent
It is not hard to see that $\ell$ is an equivalence relation on the set 
of vertices of graph $\overrightarrow{{\Bbb T}}_m$. 
\begin{dfn}\label{dfn1}
By  a cluster $C^*$-algebra ${\Bbb A}(\mathbf{x}, B)$ 
 one understands an $AF$-algebra given by
the Bratteli diagram ${\goth B}(\mathbf{x}, B)$ of the form:
\begin{equation}
{\goth B}(\mathbf{x}, B):= \overrightarrow{{\Bbb T}}_m ~\hbox{\bf mod} ~\ell. 
\end{equation}
The rank of ${\Bbb A}(\mathbf{x}, B)$  is equal to such of cluster algebra 
${\cal A}(\mathbf{x}, B)$.  
\end{dfn}
\begin{exm}\label{exm3}
\textnormal{
If  $B_T$ is  matrix (\ref{eq5}) of  Example \ref{exm2},   then ${\goth B}(\mathbf{x}, B_T)$  
 is shown Figure 3.    (We refer the reader to Section 4  for a proof.)  
Notice  that the graph ${\goth B}(\mathbf{x}, B_T)$ is  a part of 
of the Bratteli diagram of the Mundici algebra ${\goth M}$, compare 
  [Mundici 2011,  Figure 1]  \cite{Mun3}. 
}
\end{exm}
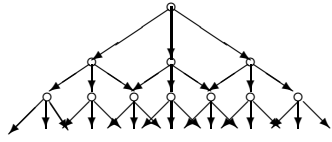
\begin{figure}[here]
\begin{picture}(100,100)(-160,130)

\put(50,200){\circle{3}}
\put(50,179){\circle{3}}
\put(20,179){\circle{3}}
\put(80,179){\circle{3}}


\put(3,166){\circle{3}}
\put(20,166){\circle{3}}
\put(36,166){\circle{3}}
\put(50,166){\circle{3}}
\put(65,166){\circle{3}}

\put(80,166){\circle{3}}
\put(98,166){\circle{3}}


\put(49,199){\vector(-3,-2){29}}
\put(51,199){\vector(3,-2){29}}
\put(50,200){\vector(0,-1){20}}


\put(19,178){\vector(-3,-2){15}}
\put(21,178){\vector(3,-2){15}}
\put(20,178){\vector(0,-1){10}}


\put(50,178){\vector(-3,-2){15}}
\put(50,178){\vector(3,-2){15}}
\put(50,178){\vector(0,-1){10}}


\put(79,178){\vector(-3,-2){15}}
\put(81,178){\vector(3,-2){15}}
\put(80,178){\vector(0,-1){10}}


\put(1.5,165){\vector(-1,-1){13}}
\put(4,164){\vector(2,-3){7}}
\put(2.5,164){\vector(0,-1){10}}

\put(19,164){\vector(-1,-1){10}}
\put(21,164){\vector(1,-1){10}}
\put(20,164){\vector(0,-1){10}}

\put(49,164){\vector(-1,-1){10}}
\put(51,164){\vector(1,-1){10}}
\put(50,164){\vector(0,-1){10}}

\put(79,164){\vector(-1,-1){10}}
\put(81,164){\vector(1,-1){10}}
\put(80,164){\vector(0,-1){10}}

\put(98,164){\vector(-1,-1){10}}
\put(100,164){\vector(1,-1){10}}
\put(98,164){\vector(0,-1){10}}


\put(65,164){\vector(-1,-1){10}}
\put(65,164){\vector(1,-1){10}}
\put(65,164){\vector(0,-1){10}}

\put(36,164){\vector(-1,-1){10}}
\put(36,164){\vector(1,-1){10}}
\put(36,164){\vector(0,-1){10}}


\end{picture}
\caption{The Bratteli diagram of  Markov's  cluster $C^*$-algebra.} 
\end{figure}

\begin{rmk}
\textnormal{
It is not hard to see that ${\goth B}(\mathbf{x}, B)$ is no longer a
tree and  ${\goth B}(\mathbf{x}, B)$  is a finite graph if and only if 
${\cal A}(\mathbf{x}, B)$ is a  finite  cluster algebra. 
}
\end{rmk}

\section{Proof}
Let $m=3(2g-2+n)$ be the rank of cluster $C^*$-algebra ${\Bbb A}(\mathbf{x}, Q_{g,n})$.
For the sake of clarity,  we shall consider   the case $m=3$ and the general  case $m\in \{ 3, 6, 9,\dots\}$
separately.

\medskip
(i)  Let    ${\Bbb A}(\mathbf{x}, B_T)$ be the cluster $C^*$-algebra of rank $3$.  
In this case $2g-2+n=1$ and either $g=0$ and $n=3$ or else  $g=n=1$.  Since
$T_{0,3}\cong \{pt\}$ is trivial,  we are left with $g=n=1$,
i.e.  the once-punctured torus $S_{1,1}$.

Repeating the argument of Example  \ref{exm2},  we get 
 the seed $(\mathbf{x}, B_T)$,   where  $\mathbf{x}=(x_1,x_2,x_3)$ and 
  the skew-symmetric  matrix $B_T$  is given by formula (\ref{eq5}).

Let us verify that matrix $B_T$ is mutation finite;  indeed,  for each $k\in \{1, 2, 3\}$
the matrix mutation  formula (\ref{eq3})  gives us $\mu_k(B_T)=-B_T$.

Therefore,  the exchange relations (\ref{eq2}) do not vary;  it is verified directly that
such relations have the form:
\begin{equation}\label{eq8}
\left\{
\begin{array}{ccc}
x_1 x_1' &=& x_2^2+x_3^2,\\
x_2 x_2' &=&  x_1^2+x_3^2,\\
x_3x_3'  &=&  x_1^2+x_2^2. 
\end{array}
\right.
\end{equation}

Consider a mutation  tree $\overrightarrow{{\Bbb T}}_3$  shown in Figure 4;
the vertices of $\overrightarrow{{\Bbb T}}_3$ correspond to the mutations
of cluster $\mathbf{x}=(x_1,x_2,x_3)$ following  the exchange rules
(\ref{eq8}).   
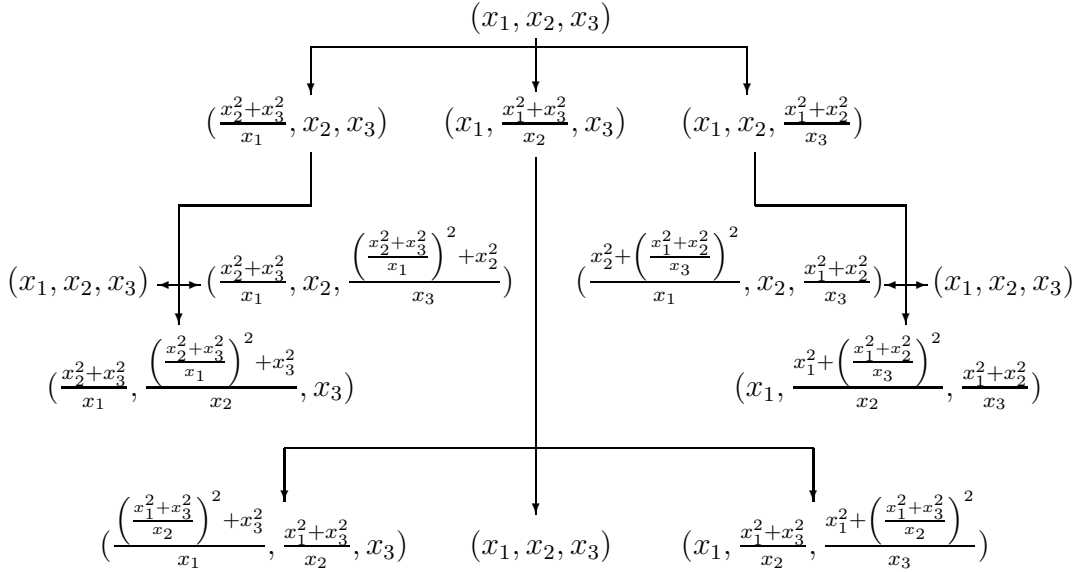
\begin{figure}[here]
\begin{picture}(300,250)(-120,70)


\put(50,290){$(x_1,x_2,x_3)$}

\put(75,285){\vector(0,-1){20}}

\put(-10,282){\vector(0,-1){18}}
\put(155,282){\vector(0,-1){18}}

\put(-10,282){\line(1,0){165}}


\put(40,250){$(x_1,{x_1^2+x_3^2\over x_2},x_3)$}
\put(-50,250){$({x_2^2+x_3^2\over x_1},x_2,x_3)$}
\put(130,250){$(x_1,x_2, {x_1^2+x_2^2\over x_3})$}


\put(-125,190){$(x_1,x_2,x_3)$}
\put(-110,150){$({x_2^2+x_3^2\over x_1},{\left({x_2^2+x_3^2\over x_1}\right)^2+x_3^2\over x_2},x_3)$}
\put(-50,190){$({x_2^2+x_3^2\over x_1},x_2, {\left({x_2^2+x_3^2\over x_1}\right)^2+x_2^2\over x_3})$}

\put(225,190){$(x_1,x_2,x_3)$}
\put(150,150){$(x_1,{x_1^2+\left({x_1^2+x_2^2\over x_3}\right)^2\over x_2}, {x_1^2+x_2^2\over x_3})$}
\put(90,190){$({x_2^2+ \left({x_1^2+x_2^2\over x_3}\right)^2\over x_1},x_2, {x_1^2+x_2^2\over x_3})$}

\put(-60,222){\line(1,0){50}}
\put(158,222){\line(1,0){57}}

\put(-10,222){\line(0,1){20}}
\put(158,222){\line(0,1){20}}

\put(215,222){\vector(0,-1){45}}
\put(-60,222){\vector(0,-1){45}}

\put(-60,192){\vector(1,0){8}}
\put(-60,192){\vector(-1,0){8}}
\put(215,192){\vector(1,0){8}}
\put(215,192){\vector(-1,0){8}}


\put(50,90){$(x_1,x_2,x_3)$}
\put(-90,90){$({\left({x_1^2+x_3^2\over x_2}\right)^2+x_3^2\over x_1},{x_1^2+x_3^2\over x_2},x_3)$}
\put(130,90){$(x_1, {x_1^2+x_3^2\over x_2}, {x_1^2+ \left({x_1^2+x_3^2\over x_2}\right)^2\over x_3})$}

\put(75,240){\vector(0,-1){135}}
\put(-20,130){\line(1,0){200}}

\put(-20,130){\vector(0,-1){20}}
\put(180,130){\vector(0,-1){20}}

\end{picture}
\caption{The mutation tree.} 
\end{figure}

The reader is encouraged to verify  that modulo a cyclic permutation of variables 
$x_1'=x_2, x_2'=x_3, x_3'=x_1$ and $x_1'=x_3, x_2'=x_1, x_3'=x_2$
one obtains (respectively) the following  equivalences of clusters:
\begin{equation}\label{eq9}
\left\{
\begin{array}{ccc}
\mu_{13}(\mathbf{x}) &=& \mu_{21}(\mathbf{x}),\\
\mu_{23}(\mathbf{x}) &=& \mu_{31}(\mathbf{x}),
\end{array}
\right.
\end{equation}
where $\mu_{ij}(\mathbf{x}):=\mu_j(\mu_i(\mathbf{x}))$;  there are no 
other cluster equivalences for the vertices of the same level of  graph  
 $\overrightarrow{{\Bbb T}}_3$.

To determine the  graph ${\goth B}(\mathbf{x}, B_T)$ one needs to take the
quotient of  $\overrightarrow{{\Bbb T}}_3$ by the  $\ell$-equivalence relations
(\ref{eq9});   since  the pattern repeats   for each level of    $\overrightarrow{{\Bbb T}}_3$,
one gets the ${\goth B}(\mathbf{x}, B_T)$ shown  in Figure 3.  
The  cluster $C^*$-algebra ${\Bbb A}(\mathbf{x}, B_T)$ is an $AF$-algebra
with the Bratteli diagram ${\goth B}(\mathbf{x}, B_T)$.

Notice that the Bratteli diagram ${\goth B}(\mathbf{x}, B_T)$ of our $AF$-algebra 
${\Bbb A}(\mathbf{x}, B_T)$ and such of the Mundici algebra ${\goth M}$ are distinct, 
compare [Mundici 2011,  Figure 1]  \cite{Mun3};  yet there is an obvious inclusion
of one diagram into another.  Namely,  if one erases a ``camel's back''  
(i.e. the two extreme sides of the diagram) in the Bratteli diagram of ${\goth M}$,
then one gets exactly the diagram in Figure 3.  Formally, if ${\cal G}$ is the Bratteli 
diagram of the Mundici algebra ${\goth M}$, 
the complement ${\cal G}-{\goth B}(\mathbf{x}, B_T)$ is a hereditary Bratteli
diagram  which  gives rise to  an ideal $I_0\subset {\goth M}$,  such that:
\begin{equation}
{\Bbb A}(\mathbf{x}, B_T)\cong {\goth M}/I_0,
\end{equation}
see [Bratteli 1972, Lemma 3.2] \cite{Bra1};    the $I_0$ is a primitive ideal {\it ibid.}, Theorem 3.8.
(It is  interesting to calculate    the group $K_0(I_0)$ in the context of  
the work of [Panti 1999]  \cite{Pan1}.)

\bigskip
On the other hand, the space $Prim ~{\goth M}$ (and hence $Prim ~{\Bbb A}(\mathbf{x}, B_T)$)
is well understood, see e.g.  [Panti 1999]  \cite{Pan1}  or   [Boca 2008,   Proposition 7]  \cite{Boc1}.  
Namely,
\begin{equation}\label{eq10}
Prim ~({\goth M}/I_0)=\{  I_{\theta} ~|~ \theta\in {\Bbb R}\},
\end{equation}
where $I_{\theta}\subset {\goth M}$ is such that ${\goth M}/I_{\theta}\cong A_{\theta}$ is
the Effros-Shen algebra [Effros \& Shen 1980]   \cite{EfSh1} if $\theta$ is an irrational number
or ${\goth M}/I_{\theta}\cong M_q$ is finite-dimensional matrix $C^*$-algebra (and an extension of such 
by the $C^*$-algebra of compact operators)  if $\theta={p\over q}$ 
is a rational number. 
(Note that the third series of primitive ideals of [Boca 2008,   Proposition 7]  \cite{Boc1}
correspond to the ideal $I_0$.) 
 Moreover,  given the Jacobson topology on $Prim ~{\goth M}$, 
there exists a homeomorphism
\begin{equation}\label{eq11}
h: Prim ~({\goth M}/I_0)\to {\Bbb R}
\end{equation}
defined by the formula $I_{\theta}\mapsto\theta$, see  [Boca 2008, Corollary 12]  \cite{Boc1}.

\bigskip  
Let $\sigma_t: {\goth M}/I_0\to {\goth M}/I_0$ be the Tomita-Takesaki
flow, i.e. a one-parameter automorphism group of ${\goth M}/I_0$,
see Section 4.  Because $I_{\theta}\subset {\goth M}/I_0$,
the image $\sigma^t(I_{\theta})$ of $I_{\theta}$ is correctly defined for
all $t\in {\Bbb R}$; the $\sigma_t(I_{\theta})$ is an ideal
of ${\goth M}/I_0$  but  not necessarily primitive.   
Since $\sigma_t$ is nothing but (an algebraic form  of)  the Teichm\"uller geodesic flow
on $T_{1,1}$  [Veech 1986]   \cite{Vee1},  one concludes that
that the family of ideals 
\begin{equation}\label{eq12}
\{\sigma_t(I_{\theta})\subset {\goth M}/I_0 ~|~  t\in {\Bbb R}, ~\theta\in {\Bbb R}\}
\end{equation}
can be taken for a coordinate system in the space $T_{1,1}\cong {\Bbb R}^2$.
In view of (\ref{eq11}) and ${\goth M}/I_0\cong {\Bbb A}(\mathbf{x}, Q_{1,1})$,   
one gets the required homeomorphism 
\begin{equation}\label{eq13}
h:  Prim~{\Bbb A}(\mathbf{x}, Q_{1,1})~\times ~{\Bbb R}\to  T_{1,1},
\end{equation}
such that the quotient algebra ${\Bbb A}(\mathbf{x}, Q_{1,1})/\sigma_t(I_{\theta})$
is  a non-commutative coordinate ring  of  the Riemann surface  $S_{1,1}$.  
\begin{rmk}\label{rmk}
\textnormal{
The family of algebras $\{{\Bbb A}(\mathbf{x}, Q_{1,1})/\sigma_t(I_{\theta}) ~|~\theta=Const, ~ t\in {\Bbb R}\}$ 
are in general pairwise non-isomorphic.  
(For otherwise all ideals $\{\sigma_t(I_{\theta}) ~|~ t\in {\Bbb R}\}$
were primitive.)  Yet   their  Grothendieck semi-groups $K_0^+$ are,  see 
[Effros \& Shen 1980]   \cite{EfSh1};    the action of $\sigma_t$ is given by the formula (see Section 4):
\begin{equation}
K_0^+({\Bbb A}(\mathbf{x}, Q_{1,1})/\sigma_t(I_{\theta}))\cong e^t({\Bbb Z}+{\Bbb Z}\theta).  
\end{equation}
}
\end{rmk}

\bigskip
(ii)  The general case $m=3k=3(2g-2+n)$ is treated likewise.  
Notice that if $d=6g-6+2n$ is  dimension of the space $T_{g,n}$,  then we have
$m-d=n$;  in particular,  rank $m$ of the cluster $C^*$-algebra  
${\Bbb A}(\mathbf{x}, Q_{g,n})$ determines completely the pair $(g,n)$ provided 
$d$ is a fixed constant.  (If $d$ is not fixed,  there is only a finite number of 
different pairs  $(g,n)$ for given rank $m$.)

Let $(\mathbf{x}, B_T)$ be the seed given by the cluster  $\mathbf{x}=(x_1,\dots, x_{3k})$ 
and the skew-symmetric matrix $B_T$.  Since matrix $B_T$ comes from a triangulation of the Riemann surface $S_{g,n}$,
$B_T$  is {\it mutation finite},   see   [Williams 2014,  p.18]  \cite{Wil1};  the  exchange relations (\ref{eq2}) 
take the form: 
\begin{equation}\label{eq14}
\left\{
\begin{array}{ccc}
x_1 x_1' &=& x_2^2+x_3^2+\dots+x_{3k}^2,\\
x_2 x_2' &=&  x_1^2+x_3^2+\dots+x_{3k}^2,\\
\vdots &&\\
x_{3k}x_{3k}'  &=&  x_1^2+x_2^2+\dots+x_{3k-1}^2. 
\end{array}
\right.
\end{equation}

One can construct the mutation tree  $\overrightarrow{{\Bbb T}}_{3k}$ using relations (\ref{eq14});
the reader is encouraged to verify,  that the   $\overrightarrow{{\Bbb T}}_{3k}$ is similar to 
the one shown in Figure 4,  except for the number of the outgoing edges at each vertex is equal to $3k$.

A tedious but straightforward calculation shows that the only equivalent clusters
at the same level of  $\overrightarrow{{\Bbb T}}_{3k}$  are the ones at the
extremities of tuples $(x_1',\dots,x_{3k}')$;  in other words, one gets the 
following system of equivalences of clusters: 
\begin{equation}\label{eq15}
\left\{
\begin{array}{ccc}
\mu_{1,3k}(\mathbf{x}) &=& \mu_{21}(\mathbf{x}),\\
\mu_{2,3k}(\mathbf{x}) &=& \mu_{31}(\mathbf{x}),\\
\vdots &&\\
\mu_{3k-1,3k}(\mathbf{x}) &=& \mu_{3k,1}(\mathbf{x}), 
\end{array}
\right.
\end{equation}
where $\mu_{ij}(\mathbf{x}):=\mu_j(\mu_i(\mathbf{x}))$.

The  graph ${\goth B}(\mathbf{x}, B_T)$ is the
quotient of  $\overrightarrow{{\Bbb T}}_{3k}$ by the $\ell$-equivalence relations
(\ref{eq15});   for $k=2$ such a graph is sketched  in Figure 5.  
The  ${\Bbb A}(\mathbf{x}, Q_{g,n})$ is an $AF$-algebra
given by  the Bratteli diagram ${\goth B}(\mathbf{x}, B_T)$.  

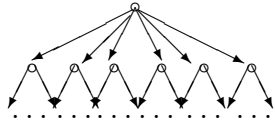
\begin{figure}[here]
\begin{picture}(100,70)(-160,130)

\put(50,200){\circle{3}}

\put(11,177){\circle{3}}
\put(27,177){\circle{3}}
\put(42,177){\circle{3}}

\put(60,177){\circle{3}}
\put(76,177){\circle{3}}
\put(94,177){\circle{3}}



\put(49,199){\vector(-2,-1){38}}
\put(48,197){\vector(-1,-1){18}}
\put(50,200){\vector(-1,-2){10}}

\put(51,199){\vector(2,-1){40}}
\put(51,199){\vector(1,-1){20}}
\put(50,200){\vector(1,-2){10}}


\put(10,178){\vector(-1,-2){8}}
\put(12,178){\vector(1,-2){8}}

\put(27,178){\vector(-1,-2){8}}
\put(28,178){\vector(1,-2){8}}

\put(42,178){\vector(-1,-2){8}}
\put(43,178){\vector(1,-2){8}}

\put(59,178){\vector(-1,-2){8}}
\put(60,178){\vector(1,-2){8}}

\put(75,178){\vector(-1,-2){8}}
\put(76,178){\vector(1,-2){8}}

\put(93,178){\vector(-1,-2){8}}
\put(94,178){\vector(1,-2){8}}


\put(88,158){$\dots$}
\put(69,158){$\dots$}
\put(51,158){$\dots$}
\put(35,158){$\dots$}
\put(20,158){$\dots$}
\put(3,158){$\dots$}

\end{picture}
\caption{The Bratteli diagram of a cluster $C^*$-algebra of rank $6$.} 
\end{figure}
\begin{lem}\label{lem1}
The set 
\begin{equation}\label{eq16}
Prim ~{\Bbb A}(\mathbf{x}, Q_{g,n})=\{ I_{\theta} ~|~ \theta\in {\Bbb R}^{6g-7+2n} ~\hbox{{\sf is generic}}\},
\end{equation}
where  ${\Bbb A}(\mathbf{x}, Q_{g,n})/I_{\theta}$ is an $AF$-algebra 
$A_{\theta}$  associated to the convergent Jacobi-Perron continued fraction of vector $(1,\theta)$,
see Section 2.3. 
\end{lem}
{\it Proof.}  We adapt the argument of    [Boca 2008,  case $k=1$]  \cite{Boc1} 
to the  case $k\ge 1$.   
Let $d=6g-6+2n$ be dimension of the space $T_{g,n}$.  Roughly speaking, the Bratteli diagram 
${\goth B}(\mathbf{x}, B_T)$ of algebra ${\Bbb A}(\mathbf{x}, Q_{g,n})$
can be cut in two disjoint pieces ${\cal G}_{\theta}$ and $ {\goth B}(\mathbf{x}, B_T)-{\cal G}_{\theta}$, 
as it is shown by   [Boca 2008,   Figure 7]   \cite{Boc1}.
The  ${\cal G}_{\theta}$ is a (finite or infinite) vertical strip of  constant ``width'' $d$,
where $d$ is equal to the  number of vertices cut from each level of  ${\goth B}(\mathbf{x}, B_T)$. 
The reader is encouraged to verify, that ${\cal G}_{\theta}$ is exactly the Bratteli diagram 
of the $AF$-algebra $A_{\theta}$ associated to the convergent Jacobi-Perron continued fraction 
of a {\it generic} vector $(1,\theta)$, see Section 2.3.  

On the other hand,  the complement  ${\goth B}(\mathbf{x}, B_T)-{\cal G}_{\theta}$ is a hereditary 
Bratteli diagram,  which defines an ideal $I_{\theta}$ of algebra  ${\Bbb A}(\mathbf{x}, Q_{g,n})$,
such that:
\begin{equation}\label{eq17}
{\Bbb A}(\mathbf{x}, Q_{g,n})/I_{\theta}=A_{\theta},
\end{equation}
see  [Bratteli 1972,  Lemma 3.2]  \cite{Bra1}.   Moreover, $I_{\theta}$ is a primitive ideal
 [Bratteli 1972,  Theorem 3.8]   \cite{Bra1}.  (An  extra care is required if $\theta=(\theta_i)$ is a rational vector;
 the complete argument can be found in [Boca 2008, pp. 980-985]  \cite{Boc1}.)   Lemma \ref{lem1} follows.

\begin{lem}\label{lem2}
The sequence of primitive ideals $I_{\theta_n}$ converges to $I_{\theta}$ 
in the Jacobson topology in $Prim~{\Bbb A}(\mathbf{x}, Q_{g,n})$
if and only if the sequence $\theta_n$ converges to  $\theta$ in the
Euclidean space ${\Bbb R}^{6g-7+2n}$.
\end{lem}
{\it Proof.}  The proof is a straightforward adaption of the argument 
in [Boca 2008, pp. 986-988]  \cite{Boc1};  we leave it as an exercise
to the reader.

\bigskip
Let $\sigma_t: {\Bbb A}(\mathbf{x}, Q_{g,n})\to {\Bbb A}(\mathbf{x}, Q_{g,n})$ be the Tomita-Takesaki
flow, i.e. the group $\{\sigma_t ~|~ t\in {\Bbb R}\}$ of  modular  automorphisms  of  algebra ${\Bbb A}(\mathbf{x}, Q_{g,n})$,
see Section 4.  Because $I_{\theta}\subset {\Bbb A}(\mathbf{x}, Q_{g,n})$,
the image $\sigma^t(I_{\theta})$ of $I_{\theta}$ is correctly defined for
all $t\in {\Bbb R}$; the $\sigma_t(I_{\theta})$ is an ideal
of ${\Bbb A}(\mathbf{x}, Q_{g,n})$ but  not necessarily a primitive ideal.   
Since $\sigma_t$ is an algebraic form  of   the Teichm\"uller geodesic flow
on the space $T_{g,n}$  [Veech 1986]   \cite{Vee1},  one concludes that
that the family of ideals: 
\begin{equation}\label{eq18}
\{\sigma_t(I_{\theta})\subset {\Bbb A}(\mathbf{x}, Q_{g,n}) ~|~  t\in {\Bbb R}, ~\theta\in {\Bbb R}^{6g-7+2n}\}
\end{equation}
can be taken for a coordinate system in the space $T_{g,n}\cong {\Bbb R}^{6g-6+2n}$.
In view of  Lemmas \ref{lem1} and  \ref{lem2},    one gets the required homeomorphism 
\begin{equation}\label{eq19}
h:  Prim~{\Bbb A}(\mathbf{x}, Q_{g,n})\times {\Bbb R}\to \{U\subseteq  T_{g,n} ~|~U~\hbox{{\sf is generic}}\},
\end{equation}
such that the quotient  algebra $A_{\theta}={\Bbb A}(\mathbf{x}, Q_{g,n})/\sigma_t(I_{\theta})$
is  a non-commutative coordinate ring  of  the Riemann surface  $S_{g,n}$.

\medskip
Theorem \ref{thm1} is proved.

\section{An analog of modular flow on  ${\Bbb A}(\mathbf{x}, Q_{g,n})$}
{\bf A.  ~Modular automorphisms $\{\sigma_t ~|~ t\in {\Bbb R}\}$.} 
Recall that the Ptolemy relations (\ref{eq4}) for the Penner coordinates $\{\lambda(\gamma_i)\}$
in the space $T_{g,n}$ are homogeneous;  in particular,  the system $\{t\lambda(\gamma_i) ~|~t\in {\Bbb R}\}$
of such coordinates will also satisfy  the Ptolemy relations.   
On the other hand,  for the cluster $C^*$-algebra ${\Bbb A}(\mathbf{x}, Q_{g,n})$
the  variables $x_i=\lambda(\gamma_i)$ and one gets an obvious isomorphism
${\Bbb A}(\mathbf{x}, Q_{g,n})\cong {\Bbb A}(t\mathbf{x}, Q_{g,n})$  for all $t\in {\Bbb R}$.
Since ${\Bbb A}(t\mathbf{x}, Q_{g,n})\subseteq {\Bbb A}(\mathbf{x}, Q_{g,n})$,  one 
obtains  a one-parmeter group of automorphisms:
\begin{equation}\label{eq20}
\sigma_t: {\Bbb A}(\mathbf{x}, Q_{g,n})\longrightarrow  {\Bbb A}(\mathbf{x}, Q_{g,n}). 
\end{equation}
By analogy with [Connes 1978] \cite{Con1}, 
we shall call $\sigma_t$ a {\it Tomita-Takesaki flow} on the cluster $C^*$-algebra 
${\Bbb A}(\mathbf{x}, Q_{g,n})$.  The reader is encouraged to verify, that $\sigma_t$ is
an algebraic form of the  {\it geodesic flow} $T^t$ on the Teichm\"uller
space $T_{g,n}$,  see [Veech 1986]   \cite{Vee1} for an introduction.
Roughly speaking,  such a flow  comes from the one-parameter group
of matrices 
\begin{equation}\label{eq21}
\left(\matrix{e^t & 0\cr 0 & e^{-t}}\right)
\end{equation}
acting on the space of holomorphic quadratic differentials
on the Riemann surface $S_{g,n}$; the latter is known to be  isomorphic 
to the Teichm\"uller space $T_{g,n}$.

\bigskip\noindent
{\bf B.  ~Connes invariant $T({\Bbb A}(\mathbf{x}, Q_{g,n}))$.}
Recall that an analogy of  the {\it Connes invariant} $T({\cal M})$ for  a $C^*$-algebra ${\cal M}$
endowed with a modular automorphism group $\sigma_t$ is the set
$T({\cal M}):=\{t\in {\Bbb R} ~|~\sigma_t ~\hbox{{\sf is inner}}\}$  [Connes 1978] \cite{Con1}.  
The group of inner automorphisms of the space $T_{g,n}$ and 
algebra ${\Bbb A}(\mathbf{x}, Q_{g,n})$ 
is  isomorphic  to  the mapping class group $Mod~S_{g,n}$  of surface $S_{g,n}$.   
The  automorphism $\phi\in Mod~S_{g,n}$ is  called {\it pseudo-Anosov},  if $\phi({\cal F}_{\mu})=\lambda_{\phi} {\cal F}_{\mu}$,
where ${\cal F}_{\mu}$ is an  invariant measured foliation and 
$\lambda_{\phi}>1$ is a constant  called {\it dilatation} of $\phi$;  the $\lambda_{\phi}$ is always
an algebraic number of the maximal degree $6g-6+2n$ [Thurston 1988]   \cite{Thu1}.  
It is known, that if $\phi\in Mod~S_{g,n}$ is pseudo-Anosov then there exists 
a trajectory ${\cal O}$ of the geodesic flow $T^t$ and a point $S_{g,n}\in T_{g,n}$,  such that the points $S_{g,n}$ and  $\phi(S_{g,n})$
belong to ${\cal O}$  [Veech 1986]   \cite{Vee1};   the ${\cal O}$ is called an {\it axis} of  the pseudo-Anosov automorphism $\phi$.   
The axis  can be used to calculate the Connes invariant 
$T({\Bbb A}(\mathbf{x}, Q_{g,n}))$ of the cluster $C^*$-algebra 
${\Bbb A}(\mathbf{x}, Q_{g,n})$;  indeed,  in view of formula (\ref{eq21})
one must solve  the following system of equations:
\begin{equation}\label{eq22}
\left\{
\begin{array}{ccc}
\sigma_t(x)  &=& e^tx \\
\phi (x)  &=&  \lambda_{\phi}x, 
\end{array}
\right.
\end{equation}
for a point $x\in {\cal O}$.  Thus $\sigma_t(x)$ coincides with the  inner automorphism 
$\phi(x)$ if and only  if $t=\log\lambda_{\phi}$.  Taking all pseudo-Anosov automorphisms  
$\phi\in Mod~S_{g,n}$,  one gets a formula for  the Connes invariant:  
\begin{equation}\label{eq23}
T({\Bbb A}(\mathbf{x}, Q_{g,n}))=\{\log\lambda_{\phi} ~|~ \phi\in Mod~S_{g,n} ~\hbox{{\sf is pseudo-Anosov}}\}. 
\end{equation}

\medskip
\begin{rmk}\label{rmk6}
\textnormal{
The Connes invariant (\ref{eq23}) says    that the family of  cluster $C^*$-algebras  ${\Bbb A}(\mathbf{x}, Q_{g,n})$
is an analog of the type $\mathbf{III}_{\lambda}$ factors of von Neumann algebras,  see  [Connes 1978] \cite{Con1}.  
}
\end{rmk}

\bigskip\noindent
{\bf Competing interests.} 
The author declares that there is no conflict of interests 
in the paper.

\bigskip\noindent
{\bf Acknowledgments.} 
It is my pleasure to thank  Ibrahim Assem and the SAG group of the Department 
of Mathematics of the University of Sherbrooke for hospitality and excellent
working conditions.  I am grateful to Ibrahim Assem, Thomas Br\"ustle,  Daniele Mundici, 
Ralf  Schiffler and Vasilisa  Shramchenko for  an introduction to  
 the wonderland  of cluster algebras and helpful correspondence;  all errors and misconceptions
in this note are solely mine. 
Finally, I thank the referee for clever comments.



\vskip1cm

\textsc{Department of Mathematics and Computer Science, St.~John's University, 8000 Utopia Parkway,  New York (Queens),  NY 11439, U.S.A.;}
~\textsc{E-mail:} {\sf igor.v.nikolaev@gmail.com}

\end{document}